\documentclass[12pt]{article}
\usepackage{graphicx}
\usepackage{amsfonts}
\newcommand{\pic}[4]{\vspace{1ex}\setlength{\unitlength}{1cm}
\begin{picture}(0,#3)(5,.5)
\put(#2){\includegraphics[#4]{#1.eps}}
\end{picture}\vspace{1ex}}

\newtheorem{prop}{Proposition}[section]

\newtheorem{theo}[prop]{Theorem}

\newtheorem{conj}[prop]{Conjecture}

\def\proof{\par\vspace{-1ex}\noindent\textit{Proof. }}
 
\def\C{\mathbb{C}}
\def\cg{\overline} 

\def\D{\mathbb{D}}
\def\DS{\displaystyle}
\def\G{\mathcal{G}}
\def\R{\mathbb{R}}
\def\S{{\mathcal{S}}}

\def\qed{\hskip1em\raise3.5pt\hbox{\framebox[2mm]{\ }}}
\def\re{\,{\rm Re}\,}
\def\im{\,{\rm Im}\,}

\newcommand{\X}[1]{X^{(#1)}}
\newcommand{\Xt}[1]{\widetilde X^{(#1)}}

\begin{document}

\begin{center} {\LARGE
    Numerical Conformal Mapping to One-Tooth Gear-Shaped Domains and Applications}
\end{center}
\begin{center}
\today
\end{center}
\begin{center}
    Philip R.  Brown\footnote{Partially supported by CONACyT grant
      166183} \\
    R. Michael Porter\footnotemark[1] 
\end{center}
 
\noindent Abstract.  We study conformal mappings from the unit disk
(or a rectangle) to one-tooth gear-shaped planar domains from the
point of view of the Schwarzian derivative, with emphasis on numerical
considerations. Applications are given to evaluation of a singular
integral, mapping to the complement of an annular rectangle, and
symmetric multitooth domains.


\medskip

\section{Introduction \label{sec:intro}}

In \cite{BrP2} we initiated a study of conformal mappings to a
\textit{gearlike} domain with a single tooth: a starlike open
set in the complex plane bounded by arcs of two circles centered at
the origin and segments of two lines passing through the origin.
The approach was to examine the Schwarzian derivative of such a
mapping (a rational function of a particular form).  Relationships
were found between the auxiliary parameters in the Schwarzian derivative
and the geometry of the gearlike domain, and the set of parameters producing
a univalent mapping (``region of gearlikeness'') was determined.

Using the theoretical results in \cite{BrP2}, we focus here on the
computational aspects of the conformal mappings to one-tooth gear
domains from both a disk and a rectangle. We begin in
section~\ref{sec:gearlike} with the relationship between the
parameters in the Schwarzian derivatives for these two cases.  In
section~\ref{sec:formulas} we work out the transformations or
normalizations necessary for the solution of the Schwarzian equation
to produce a gear (rather than just a ``pregear'' in the terminology
of \cite{BrP2}).  In section~\ref{sec:compu} we describe a variety of
approaches for carrying out the computations, and report some
numerical results in section~\ref{sec:numerical}. This includes a
``map'' of the internal structure of the region of gearlikeness.

Finally, in Section~\ref{sec:appl} we give three applications of our
results on gear mappings: evaluation of a singular integral, mapping
to the complement of an annular rectangle, and symmetric multitooth 
domains.

An Appendix contains some background material on one of the numerical
methods employed.

\section{One-tooth gear domains\label{sec:gearlike}}

General facts stated below concerning conformal mapping, Schwarzian
derivatives, etc., may be found in many standard texts such as
\cite{DT,Hen,Hi,Neh}.  A \emph{one-tooth gear domain} is a topological
quadrilateral in the complex plane which is the union of an open disk
and a sector of a concentric disk of larger radius.  Thus up to an
affine equivalence $az+b$, a one-tooth gear domain is the standard
domain $G_{\beta,\gamma}$ with vertices at $w_1=\beta e^{i\gamma}$,
$w_2=e^{i \gamma}$, $w_3=e^{-i \gamma}$, $w_4= \beta e^{-i \gamma}$,
where the edges $[w_1,w_2]$ and $[w_3,w_4]$ are straight segments
called the \emph{tooth edges} of the gear, and the edges from $w_2$ to
$w_3$ and from $w_4$ to $w_1$ are arcs of the circles $\{|w|=1\}$ and
$\{|w|=\beta\}$ subtending angles in the ranges $\gamma<\arg
w<2\pi-\gamma$ and $-\gamma<\arg w<\gamma$ respectively.  We say that
$\beta$ is the \emph{gear ratio} and $\gamma$ is the \emph{gear
  angle}.

From our perspective, the first thing one needs is to have available
is the Schwarzian derivative $S_f=(f''/f')'-(1/2)(f''/f')^2$ of the
conformal mapping.  When $S_f$ is known, a unique solution $f$ is
determined by the triple $J_f(0)$, where the \emph{2-jet} of $f$ at
any point $z$ is defined to be
\begin{equation}  \label{eq:def2jet}
   J_f(z) = (f(z),f'(z),f''(z)),
\end{equation}
as long as we assure that $f'(0)\not=0$.

As we note in Proposition 2.1 of \cite{BrP}, if $f(0)$ is the tooth
center (the common center of the concentric circles containing the
circular boundary arcs of the gear) and $f'(0)=1$, we have
\begin{equation}  \label{eq:standard2jet}
   J_f(0) = (0,1,2(\cos t_2-\cos t_1)) 
\end{equation}
where $0<t_1<t_2<\pi$ and $e^{\pm i\pi t_1}$ and $e^{\pm i\pi t_2}$
are the prevertices. However, when we solve for $f$ given $S_f$ and
with the normalization (\ref{eq:standard2jet}), $f(0)$ does not
generally turn out to be the tooth center.  Unless we know precisely
what 2-jet to use, the solution $f$ will be a M\"obius transformation
of a one-tooth gear domain, which we term a \emph{pregear} (see Figure
\ref{fig:pregearillustration}).  Although pregears are less rigid
objects than gears, they have some very restrictive properties. Note
that the tooth edges of any pregear may be uniquely identified by having 
different interior angles at their endpoints. Applying Euclidean transformations, we
restrict the discussion to pregears which are symmetric in $\R$ and
have no vertices on $\R$. In \cite{BrP} the following was shown.
 
\begin{prop} \label{prop:pregearcondition} (a) Let $D$ be a circular
  quadrilateral with the above symmetries, having two interior angles
  equal to $\pi/2$ and two interior angles equal to $3\pi/2$.  Assume
  that one tooth edge of $D$ lies in the upper and the other in the
  lower half-plane.  Then $D$ is a pregear if and only if the full
  circles $C^+$, $C^-$ containing the tooth edges intersect in two
  points.

  (b) Let $D$ be a pregear.  Then $D$ is a gear if and only if its
  tooth edges are straight, or equivalently, if the non-tooth edges
  are arcs of concentric circles.
\end{prop} 

\subsection{Conformal mapping from disk to a one-tooth gear\label{subsec:diskschw}}

Let $\D=\{z\in\C\colon\ |z|<1\}$ denote the unit disk, and let
$f\colon\D\to G_{\beta,\gamma}$ be a conformal mapping. Suppose the
prevertices $z_i=f^{-1}(w_i)$ are located at points
of the form $e^{\pm it_1}$ and $e^{\pm it_2}$, for
$0<t_1<t_2<\pi$. The expression for the Schwarzian derivative $\S_f$
of $f$ as a rational function $S_f=R_{t_1,t_2,\lambda}$ in terms of
$t_1$, $t_1$, and an auxiliary parameter $\lambda$ was worked out
explicitly in \cite{BrP2}. There it was also noted that by means of
precomposition of $f$ with a M\"obius transformation 
\begin{equation}  \label{eq:Tq}
  T_q(z)= \frac{z-q}{1-qz}
\end{equation} 
($-1<q<1$) which leaves $\D$ invariant, the prevertices can be
symmetrized with respect to the imaginary axis. Then (writing now $f$ in
place of $f\circ T_q$) the Schwarzian derivative takes the form
\begin{equation}\label{eq:SfRtlambda0}
S_f=R_{t,\lambda}
\end{equation}
where
\begin{eqnarray}  \label{eq:Rtlambda}
    \frac{1}{2}R_{t,\lambda}(z) &=& \psi_{0,t}(z) - \lambda
    \psi_{1,t}(z)
\end{eqnarray}
and
 \begin{eqnarray}  \label{eq:psi0tpsi1t} 
  \psi_{0,t}(z) 
   &=&    \frac{ (\sin^2t)(z^4-(16\cos t)z^3 + (4+2\cos2t) z^2 - (16\cos t) z + 1) }
           {2(z^4-(2 \cos2t)z^2+1)^2}, \nonumber \\
  \psi_{1,t}(z)  &=&    \frac{-8\cos t}{z^4-(2\cos2t)z^2+1}.  \label{eq:psit}
\end{eqnarray}
One advantage of the symmetrized form is that the conformal $M(t)$ of
the gear domain is easily related to the single parameter $t$ via an
elliptic integral, as we will see in a moment. In \cite[section
3]{BrP2} this was exploited to derive qualitative relationships
between the pair of geometric parameters $\beta,\gamma$ that prescribe
the gear and the mapping parameters $t,\lambda$. In this paper we are
more concerned with quantitative information concerning this
relationship.

\subsection{Conformal mapping from rectangle to gear}\label{subsec:rectangle}
 
In the study of circular quadrilaterals with two symmetries in
\cite{BrP}, mappings from a rectangle to the quadrilateral were
investigated as well as mappings from a disk. For numerical work such
mappings have certain advantages over mapping from the disk, among
them the fact that certain Schwarzian derivatives are real on the
boundary (or on horizontal or vertical sections of the rectangle).
The relationships between conformal mappings to a gear from a disk and
from a rectangle given below are somewhat more complicated than in the
situation discussed in \cite{BrP} but follow the same general line of
reasoning.

As in \cite{BrP}, we let $E$ denote the elliptic integral
\begin{equation}  \label{eq:ellipticintegral}
 E(z) = \int_0^z \frac{dz}{ \sqrt{
         (z-e^{it})(z+e^{-it})(z+e^{it})(z-e^{-it})} },\quad |z|<1.
\end{equation}
This is a particular case of a Schwarz-Christoffel mapping \cite{DT},
and the image
$R=E(\D)=[-\omega_1,\omega_1]\times[-\im\omega_2,\im\omega_2]$ is a
rectangle centered at the origin.  Here $\omega_1$, $\omega_2$ are
half-periods of the related Weierstass $\wp$-function, with
$\omega_1>0$, $\im \omega_2=|\omega_2|$, and $\tau=\omega_2/\omega_1$.
Thus the composition $f\circ E^{-1}$ is a vertex-preserving conformal
mapping from a rectangle to the gear $G_{\beta,\gamma}$, whose
conformal module \cite{LV} is by definition the ratio
\begin{equation}  \label{eq:M(t)}
  M(t) = \frac{ \im E(i) }{ E(1) }.
\end{equation}
of the sides of the rectangle.
 
\begin{figure}[!b]
  \centering
   \pic{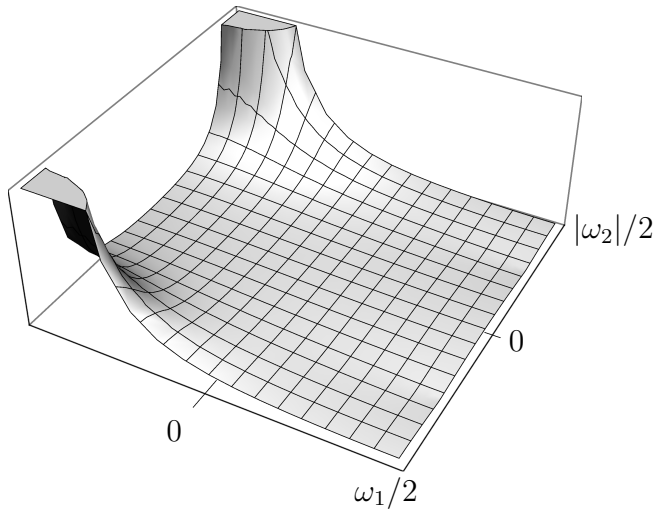}{1,0}{6.5}{ }
  \caption{  $|\varphi_{\tau,\mu}|$ for Schwarzian derivative  
(\ref{eq:varphi})}
  \label{fig:rectangleschwarzian}
\end{figure}

As customary, write $e_i=\wp(\omega_i)$, $i=1,2,3$ where
$\omega_3=\omega_1+\omega_2$. These values are real and satisfy $e_2
<e_3<e_1$.  Also as in \cite{BrP} for definiteness we will have
$\omega_1$, $\omega_2$ normalized so that $e_1-e_2=4$, so the periods
are in fact uniquely determined by their ratio $\tau$. The change of
variable
\begin{equation}  \label{eq:zetafromz}  
  \zeta= \frac{E(z)}{2}
\end{equation}
maps $\D$ onto the subrectangle
\begin{equation}  
  R_0  = \{\zeta\colon\ -\frac{\omega_1}{2} < \re \zeta < \frac{\omega_1}{2} ,\
           -\im\frac{\omega_2}{2}  < \im \zeta < \im\frac{\omega_2}{2}  \} 
\end{equation}
of $R$, and $\pm e^{\pm it}$ are mapped to the vertices of $R_0$. We define
\begin{equation}  \label{eq:varphi}
   \varphi_{\tau,\mu}(\zeta) = 
    - 4\left(\wp\bigg(\zeta+\frac{\omega_1+\omega_2}{2}\bigg)
    + \wp\bigg(\zeta+\frac{\omega_1-\omega_2}{2}\bigg)\right) + 4\mu
\end{equation}
where $\mu$ is real (Figure \ref{fig:rectangleschwarzian}).  
It is straightforward to verify, using standard properties of
Weierstrass elliptic functions \cite{WW} that $\varphi_{\tau,\mu}(\zeta)$ is
real when $\zeta\in\R$ and also when $\zeta\in\partial R_0$ (although not for
imaginary $\zeta$ in general).  The composition
\begin{equation}  \label{eq:def:g}
  g(\zeta) = f\circ E^{-1}(2\zeta)
\end{equation}
is our mapping $g\colon R_0\to G_{\beta,\gamma}$.
 
\begin{prop}\label{prop:recgear}
  Let $G=G_{\beta,\gamma}$ be a gear domain and let $g\colon R_0\to G$ be a conformal
  mapping which respects vertices, taking the horizontal edges of $R_0$
  to the tooth edges of $G$ and the left vertical edge of $R_0$ to the
  inner circumference of $G$. Then the Schwarzian derivative of $g$ is
  equal to $\S_g=\varphi_{\tau,\mu}$ for some values of $\tau$ and
  $\mu$.
\end{prop}

\proof The argument is quite similar to that of \cite[eq.\ (12)]{BrP},
so we will only mention the salient points.  One extends $g$ by
reflection across the edges of $R_0$ and notes that the composition of
reflections along opposite edges in the $\zeta$-plane is a translation
of magnitude $2\omega_1$ or $2\omega_2$. By the Chain Rule for
Schwarzian derivatives, one sees that $\S_g(\zeta)$ is a function which has
periods of both of these magnitudes.  Due to the fact that $g$ sends
the right angles at the vertices $\omega_3/2=(\omega_1+\omega_2)/2$
and $\cg{\omega_3}$ to the right interior angles of $G$, a calculation
of the series expansion centered at these vertices shows that $g$ is
holomorphic at these vertices.  In contrast, due to the angles of
$3\pi/2$ at the remaining vertices of $G$, the Schwarzian derivative
$S_g$ has double poles at the left vertices of $R_0$, with
singularities
\[  \frac{-4}{(\zeta+\omega_3/2)^2},\ 
    \frac{-4}{(\zeta+\cg{\omega_3}/2)^2}.
\]
One checks that the function $\varphi_{\tau,0}$ given explicitly by
(\ref{eq:varphi}) is an elliptic function with the same singularities,
so the difference is a constant, which is then seen to be real due to the
symmetry. \qed

The parameter $\mu$ plays a role similar to that of the parameter of
the same name in \cite{BrP}, just as does $\lambda$ for our mappings of
the disk, although the Schwarzians are now different functions in the
present work.

\subsection{Relation between disk and rectangle Schwarzians}
 
\begin{prop}
  Let $0<t<\pi/2$, let $\tau/i>0$, and let $f\colon\D\to\C$ and
  $g\colon R_0\to\C$ be conformal mappings with Schwarzian derivatives
\[ S_f(z) = R_{t,\lambda}(z), \quad S_g(\zeta) = \varphi_{\tau,\mu}(\zeta)
\]
respectively.  Suppose that $g(\zeta)=f(z)$ where (\ref{eq:zetafromz}) holds.
Then $\tau/i=M(t)$ and 
 \begin{equation}   \label{eq:mufromlambda}
   \mu = 16 \lambda \cos t + \frac{3 + \cos2t}{6} .
\end{equation}
\end{prop}

\proof The relation between $t$ and $\tau$ is evident by the definition of conformal
module. Define
\[ \varphi_1(z) = \wp (\zeta+\frac{\omega_1+\omega_2}{2}), \quad
   \varphi_2(z) = \wp (\zeta+\frac{\omega_1-\omega_2}{2}),
\]
so $\varphi_{\tau,0}(\zeta) = -4(\varphi_1(z)+\varphi_2(z))$.  From
the symmetries and the well-known fact \cite{WW} that $\wp\colon
(\C\bmod\{2\omega_1,2\omega_2\})\to\C\cup\{\infty\}$ is a 2-to-1
covering branched at the vertices of $R$, we see that $\varphi_1$ and
$\varphi_2$ map $\D$ onto the upper and lower half planes,
respectively. By the Reflection Principle, it can be deduced that
$\varphi_1$ and $\varphi_2$ are M\"obius transformations. Further,
they send the quadruple $(e^{it},-e^{-it},-e^{it},e^{-it})$ to
\[ (e_3,\ e_2,\ \infty,\ e_1), \quad (e_2,\ e_3,\ e_1,\  \infty)
\]
respectively. Note that cross ratios of the above quartets are equal
due to the relation $e_1+e_2+e_3=0$. From this it follows that
$\varphi_2=(e_2\varphi_1+e^2+e_1e_3)/(\varphi_1-e_2)$, and then a
straightforward calculation verifies that
\begin{eqnarray*}
  \varphi_1(z) &=& e_3  + i(e_1-e_3)(\cot t)\frac{z-e^{it}}{z+e^{it}},\\
  \varphi_2(z) &=& e_3  - i(e_1-e_3)(\cot t)\frac{z-e^{-it}}{z+e^{-it}}.
\end{eqnarray*}
We now substitute these formulas in (\ref{eq:varphi}) to obtain a formula
for $\varphi_{\tau,\mu}(\zeta)$ as a function of $z$,
\begin{equation}
 \varphi_{\tau,\mu}(\zeta) = \frac{ (\frac{8}{3}\cos2t)z^2+
      (\frac{40}{3}\cos2t\cos t-8\cos t)z + \frac{8}{3}\cos2t}{z^2+ (2\cos2t)z +1 } + 4\mu.
\end{equation}
 The Chain Rule for Schwarzian derivatives says that
\[  S_f(z) = S_g(\zeta)\left(\frac{E'(z)}{2}\right)^2 + S_E(z)
\]
and the terms involving $E$ are easily evaluated from the definition
(\ref{eq:ellipticintegral}) as in \cite{BrP}.  After simplification, we find that
\[ S_f(z) = 
  \frac{  c_0 y  +c_1 z  +   c_2 z^2+c_3 z^3   + c_4z^4  }{
3 \Delta^2}
+\frac{\mu}{\Delta}
\]
where $c_0=c_4=-2\cos2t$, $c_1=c_3=12(\cos3t-\cos t)$, $c_2=5-\cos4t$,
$\Delta=z^4-(2\cos2t)z^2+1$. On comparison with
(\ref{eq:Rtlambda})--(\ref{eq:psit}), we
arrive at the desired relationship between $\lambda$ and $\mu$.  \qed
 
\section{Formulas for calculating a gear mapping\label{sec:formulas}}

Our numerical work will rest heavily  on the following fact.

\begin{figure}[!b]   \centering
  \pic{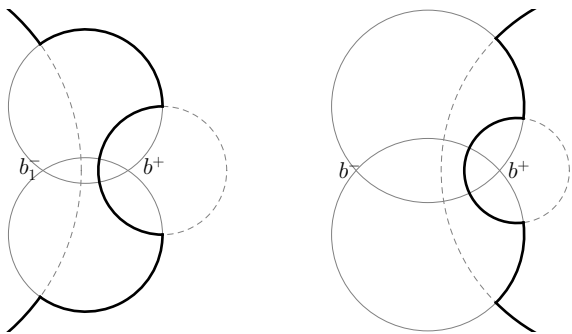}{1,0}{4}{scale=.5}
  \caption{Structure of pregears as in the proof of Proposition
    \ref{prop:pregeartogear}.  Note that precisely one of $b^{\pm}$ is
    interior to the pregear.}
\label{fig:pregearillustration}
\end{figure}
 
\begin{prop}\label{prop:pregeartogear}
  Let $D$ be a pregear symmetric in $\R$.  Suppose that $D$ is not a
  gear. Let $p_{-1},p_1\in\partial D$ ($p_{-1}<p_1$) be the extreme
  points of the real interval $D\cap\R$ (i.e., the midpoints of
  the non-tooth edges of $D$). Let $\rho$ be the common radius of the 
  tooth edges, and let $p\in\partial D$ be an interior point of
  either tooth edge.  Let $v$ be the inward pointing unit normal to the tooth
  edge at $p$. Define
\[ c = p+rv  ,\quad d= \left( \rho^2 - (\im   c)^2\right)^{1/2},\quad b^\pm=\re c\pm d,
\]
and 
\[  T_{b^\pm}(z) = \left\{\begin{array}{rl} 
    \DS -\frac{z-b^-}{z-b^+}  & 
 \mbox{ if } p_{-1}<b^-<p_1,\\[3ex]
   \DS \frac{z-b^+}{z-b^-}  & 
 \mbox{ if } p_{-1}<b^+<p_1.
 \end{array}\right.
 \]
 Then the image $G=T_{b^\pm}(D)$ is a gear with gear center at the
 origin. Its gear parameters are given by
\begin{eqnarray*}
  \beta &=&   \frac{T_{b^\pm}(p_1)}{T_{b^\pm}(p_{-1})} ,\\
  \gamma &=& \arg T_{b^\pm}(p).
\end{eqnarray*}
\end{prop} 
\proof By construction, $c$ is the center of a tooth edge.  By
Proposition \ref{prop:pregearcondition}, the tooth edges must meet,
which is impossible unless $|\im c| < \rho$ (Figure
\ref{fig:pregearillustration}). Thus $d$ is real (we assume positive),
and the intersection points of the tooth edges are thus $b^{\pm}$, with
$b^-<b^+$. The
definition of $T_{b^\pm}(z)$ is arranged so that $T_{b^\pm}$ leaves
$\R\cup\{\infty\}$ invariant, it has its pole exterior to $D$, and
$T_{b^\pm}'(x)>0$ for $x\in\R$.  The points $b^\pm$ are sent by
$T_{b^\pm}$ to $0,\infty$, and hence the images of the edges are straight, so $G$
is a gear by Proposition \ref{prop:pregearcondition}, which further
implies that the images of the remaining edges are circles centered at
the origin of radii $-T_{b^\pm}(p_{-1})$ and $T_{b^\pm}(p_1)$.  Since
$T_{b^\pm}(p)$ lies on the upper tooth edge and this edge prolongs to
pass through the origin, the formulas stated for $\beta$, $\gamma$ hold. \qed
    
An operation equivalent to changing the 2-jet at the origin may be
effected with the aid of a self-mapping $T_q$ of $\D$ defined by
(\ref{eq:Tq}). We use the following, which is verified by a direct
application of the Chain Rule.

\begin{prop}\label{prop:2jetcomp} 
 Let $f$ and $g$ be two mappings such that $J_f(z_0)=(a_0,a_1,a_2)$ and
 $J_g(a_0)=(b_0,b_1,b_2)$. Then
\[ J_{g\circ f}(z_0) = (b_0, a_1b_1,\ a_1^2b_2+a_2b_1).
\]
In particular,
\[ J_{f\circ T_q}(0)=
\left(f(-q),\ (1\!-\!q^2)f'(-q),\ (1\!-\!q^2)\left((1\!-\!q^2)f''(-q)+2q\,f'(-q)\right)\right).
\]
\end{prop}

The mapping $f\circ T_q$ does not in general fix the origin. To
achieve this condition we must use $h=T_{b^\pm}\circ f\circ T_q$ where
$T_{b^\pm}$ was defined in Proposition \ref{prop:pregeartogear}. If
$p_{-1}<b^-<p_1$, then the 2-jet of the M\"obius transformation
$T_{b^\pm}$ at $b^-$ is
\[  J_{T_{b^\pm}}(b^-)= \left(0,\ \frac{1}{b^+-b^-},\ \frac{-2}{(b^+-b^-)^2} \right).
\]
Otherwise, if $p_{-1}<b^+<p_1$, then the 2-jet of the M\"obius transformation
$T_{b^\pm}$ at $b^+$ is
\[  J_{T_{b^\pm}}(b^+)= \left(0,\ \frac{1}{b^+-b^-},\ \frac{2}{(b^+-b^-)^2} \right).
\]
We now apply Proposition \ref{prop:2jetcomp} to obtain
\begin{eqnarray} \label{eq:2jetzerotozero}
 J_h(0) &=& \left(0,\ \frac{1-q^2}{b^+-b^-}f'(-q),\  
  \frac{\pm 2(1-q^2)^2}{(b^+-b^-)^2}f'(-q)^2 +   \right. \nonumber\\
 &&  \quad\quad\quad  \left. \frac{1-q^2}{b^+-b^-}((1-q^2)f''(-q)+2qf'(-q) \right) 
\end{eqnarray}
and as a result we have the following full description of the
normalized gear mapping.
\begin{prop} \label{prop:zerotozero} Given $t,\lambda$, let $f$ denote
  the solution of $\S_f=R_{t,\lambda}$ normalized by $J_f(0)=(0,1,0)$.
  Suppose that the image $f(\D)$ is a pregear.  Let $q=-f^{-1}(b^-)$
  or $q=-f^{-1}(b^+)$
  where the pregear center $b^-\in\R$ or $b^+\in\R$ is as described in Proposition
  \ref{prop:pregeartogear}.  Then the solution $h$ of the Schwarzian
  differential equation $\S_h=(R_{t,\lambda}\circ T_q)(T_q')^2$
  normalized by the 2-jet of (\ref{eq:2jetzerotozero}) is a gear (not
  only a pregear) mapping satisfying $h(0)=0$, $h'(0)>0$.
\end{prop}

In order to apply Proposition \ref{prop:zerotozero}, it is necessary to know the
radius $\rho$ in Proposition \ref{prop:pregeartogear}. The following result
will fill this need.

\begin{prop} \label{prop:curv}
 Let $f$ be holomorphic near $z_0=e^{it_0}$, $f'(z_0)\not=0$, and suppose that for $|z|=1$
near $e^{it_0}$,
$f(z)$ lies in an arc of a circle of some radius $\rho$. Then the curvature of this arc is
\[  \frac{1}{\rho} = \frac{\pm 1}{|f'(z_0)|}\re\left(1+z_0\frac{f''(z_0)}{f'(z_0)} \right).
\]   
\end{prop} 
 A simpler version of this result was used in \cite{BrP,KP2}.  The choice of sign in this
curvature formula corresponds to whether $f(rz_0)$ enters or leaves the circle containing the
arc when $r<1$ increases to $r>1$.

\proof It is sufficient to consider the special case that $z_0=1$ (and then apply the result
to $f(e^{it_0}z)$).
 Let $c$ denote the center of the circle $|w-c|=\rho$ containing the image arc.  By the
Reflection Principle,
\[ (f(x)-c)(\overline{f(1/x)-c}) = \rho^2
\]
for $1-\epsilon<x<1+\epsilon$. From the derivative
 \[ f'(x)(\overline{f(1/x)-c})  - {1
   \over x^2}(f(x)-c)\overline{f'(1/x)} = 0
\]
we see that $f'(1)(\overline{f(1)-c})=(f(1)-c)\overline{f'(1)}$, or equivalently
\begin{equation}\label{imarg}
 \im \overline{f'(1)}(f(1)-c)=0,\quad     \arg(f(1)-c) \equiv \arg f'(1) \bmod \pi.
\end{equation}
Take a further derivative and evaluate at $x=1$:
\[ f''(1)(\overline{f(1)-c}) 
   - 2|f'(1)|^2  + 2\overline{f'(1)}(f(1)-c)+\overline{f''(1)}(f(1)-c) =0,
\]
from which it follows that
\[  2\re \overline{f''(1)}(f(1)-c) - 2|f'(1)|^2 + 2\re \overline{f'(1)}(f(1)-c) = 0
\]
so
\[  \re( (f''(1)+f'(1)) \overline{(f(1)-c)} ) = |f'(1)|^2.
\]
As a consequence of (\ref{imarg}),
\begin{eqnarray*}
  \re (f''(1)+f'(1))\overline{(f(1)-c)} &=& \re\left(
   (f''(1)+f'(1))|f(1)-c|\frac{\pm|f'(1)|}{f'(1)}\right) \\
 &=& \pm |f(1)-c| \re\left( (f''(1)+f'(1))\frac{|f'(1)|}{f'(1)} \right).
\end{eqnarray*}
We therefore have 
\[ \frac{\pm1}{ |f(1)-c|} = \frac{\re\left((f''(1)+f'(1))\frac{|f'(1)|}{f'(1)}\right) } 
     {|f'(1)|^2} = \frac{1}{|f'(1)|} \re\left( 1+\frac{f''(1)}{f'(1)}\right).
\qed
\]
 
\section{Numerical computation of gear parameters\label{sec:compu}}

\subsection{Auxiliary second-order linear ODE}

The basic facts relating the second-order linear differential equation
$2y''+R_{t,\lambda}y=0$ to the Schwarzian derivative of a conformal
mapping are widely known \cite{Neh} and have been applied in many
contexts.  A conformal $f$ mapping with Schwarzian derivative
$R_{t,\lambda}$ is obtained as a quotient
\begin{equation}  \label{eq:y2/y1}
  f = \frac{y_2}{y_1}
\end{equation}
 of any two
linearly independent solutions of the differential
equation
\begin{equation}  \label{eq:SL}
  2y'' +  R_{t,\lambda}y =0
\end{equation}
or equivalently,
as an antiderivative of $y_1^{-2}$ (we assume $y_1$ nonvanishing).

Thus given parameters $t,\lambda$, it is a straightforward matter to
compute $f$ numerically such that
\begin{equation}    \label{eq:SfRtlambda}
   \S_f = R_{t,\lambda}.
\end{equation}
A particular solution will depend on the normalization chosen for
$y_1$, $y_2$ at a given base point. In general the image $D=f(\D)$
will not be a gear.  Assuming it is at least a pregear, in order to
apply Proposition \ref{prop:curv} and Proposition \ref{prop:pregeartogear},
we must be prepared to calculate the 2-jet $J_f$ at certain boundary
points, which may be accomplished by the relations
\begin{equation}   \label{eq:2jetfromode}
    f=\frac{y_2}{y_1},\quad  f'=  \frac{1}{y_1^2} \quad
 f'' = \frac{-2y_1'}{y_1^3} 
\end{equation}
where we assume that the constant $y_1y_2'-y_2y_1'$ is equal to 1.
 
\subsection{Computational procedures for the disk\label{subsec:compdisk}}

We describe here some aspects of the calculation for mappings defined
in the disk $\D$.  Similar considerations apply to mappings defined
in the rectangle $R_0$, as well as further observations we will make
in \ref{subsec:comprec}. First we discuss calculation along radii.

Consider a fixed value of $t_0$.  When we write equation (\ref{eq:SL})
by parametrizing along the radius from $0$ to $z_0=e^{it_0}$, it takes
the form
\begin{equation}  \label{eq:SLt}
 \eta''(r) + e^{2it_0}  R_{t,\lambda}(re^{it_0})\eta(r) =  0 
\end{equation}
with $\eta(r)=y(re^{it_0})$.
We use $t_0=0,\pi/2,\pi$, i.e., we calculate the values of $f$ along
the rays from 0 to $1$, $i$, $-1$.
We consider two general approaches to the computation.

\subsubsection{Calculation via the 2-jet.\label{subsubsec:diskradial}} 
The first, direct approach is simply to solve the equations
(\ref{eq:SLt}) numerically along radii. The values of the solutions at
$r=1$ together with relation (\ref{eq:2jetfromode}) and Proposition
\ref{prop:curv} and Proposition \ref{prop:pregeartogear} permit us to
calculate the transformation $T(z)$ which sends the pregear to a gear.
Then $F=T\circ f$ is a gear mapping with $S_f=R_{t,\lambda}$.
 
\subsubsection{Calculation via the spectral parameter power series.\label{subsubsec:diskspps}} 
In the second approach, using (\ref{eq:Rtlambda})--(\ref{eq:psit}) we
rewrite (\ref{eq:SL}) as
\begin{equation}  \label{eq:SLpsi}
 y'' + \psi_{0,t}y = \lambda  \psi_{1,t}y
\end{equation}
 and again integrate along the three rays from 0 to $1$, $i$, $-1$.  We
are following very closely the notation used in \cite{BrP,KP2}, but
now the functions $\psi_0,\psi_1$ are somewhat different since we are
considering a different class of quadrilaterals. Then we make use of
the SPPS representation of the solutions of (\ref{eq:SLpsi}), as was
done in \cite{BrP,KP2}.  A summary of the main result establishing
this representation is given in the Appendix. The procedure for gear
mappings is as follows.
 
\begin{itemize}
\item Calculate the SPPS integrals routinely $\widetilde X^{(n)}$,
  $X^{(n)}$ in terms of the data $\psi_{0,t}$, $\psi_{1,t}$. These are
  functions on $0\le r\le 1$ which depend on $t$ but not on
  $\lambda$. These indefinite integrals produce the coefficients of
  power series (\ref{eq:SPPSseries}) in $\lambda$ defining linearly
  independent solutions $\eta_1,\eta_2$ of (\ref{eq:SLpsi}) in $\D$,
  normalized by 1-jets $(\eta_1(0),\eta_1'(0))$ and
  $(\eta_2(0),\eta_2'(0))$ equal to $(1,0)$ and $(0,e^{it_0})$,
  respectively.  This normalization corresponds via (\ref{eq:y2/y1})
  to the conformal mapping $f$ normalized by the 2-jet
\[   J_f(0) = (0,1,0)
\]
at $z=0$.
\item Evaluate the coefficients $\Xt{n}$, $\X{n}$ of the power series at $r=1$. This
  defines $\eta_1(1)$, $\eta_2(1)$ as functions of $\lambda$ (separately for each 
  of the three rays).
\item Express the 2-jet $J_f$ as a function of $\lambda$ by means of (\ref{eq:2jetfromode}). 
\end{itemize}
  
 Some variants of these approaches will be mentioned later.

\subsection{Computational procedure for the rectangle \label{subsec:comprec}}

\subsubsection{Condition to make pregear into a gear}

Analogously to the mappings of the disk, the $\S_g$ can be expressed in
terms of solutions of the ordinary differential equation
\begin{equation}  \label{eq:odeR0}
  2y''(\zeta) + \varphi_{\tau,\mu}(\zeta) y(\zeta) = 0.
\end{equation}  
Two particular solutions are normalized by $J_{y_1}(0)=(1,0)$,
$J_{y_2}(0)=(0,1)$.  The quotient $g_0=y_2/y_1$, normalized by
$J_{g_0}(0)=(0,1,0)$, generally does not have a gear as image; the gear mapping
must be sought within the general solution of $S_g=\varphi_{\tau,\mu}$, i.e.
\begin{equation}
  \label{eq:gdef}
  g = \frac{ a\frac{y_2}{y_1}+b}{c\frac{y_2}{y_1}+d} 
   = \frac{by_1+ay_2}{dy_1+cy_2}
\end{equation}
with $ad-bc=1$.  Thus
\[ g(0)=\frac{b}{c},\quad g'(0)=\frac{ad-bc}{d^2},\quad
 g''(0)=\frac{-2(ad-bc)c}{d^3}
\]
We may assume $a,b,c,d\in\R$ because we are only interested in the
case that $g'(\zeta)\in\R$ for $\zeta\in\R$. We may further suppose
that $g(0)=0$, $g'(0)$=1 without affecting whether the image is a gear
or not, so with this simplification we have $b=0$ and $a=d$. Therefore we write
\begin{equation}   \label{eq:ggeneral}
  g =  \frac{y_2}{y_1+\alpha y_2}
\end{equation}
where $\alpha=c/a\in\R$, and the question can be restated as how to
choose $\alpha$ to make the pregear a gear.

\begin{prop}  \label{prop:gearconditiong}
  Let the mapping $g$ be determined by (\ref{eq:ggeneral}) where
  $y_1$, $y_2$ are normalized solutions of (\ref{eq:odeR0}) in $R_0$.
  Suppose that the image $g_0(R_0)$ of $g_0=y_2/y_1$ is a pregear.
  Then $g$ is a gear mapping if and only if\/ $\alpha$ is a root of
  the quadratic equation
\[ \left( \im y_2\cg{y_2'}\,|_{\omega_3/2}\right)\alpha^2  +
   \left( \im(y_1\cg{y_2'}+ y_2\cg{y_1'})\,|_{\omega_3/2}\right)\alpha   +
   \im y_1\cg{y_1'}\,|_{\omega_3/2}  = 0
\]
and the image $g(\D)$ is bounded. If the quadratic equation has
no roots, then $g_0$ does not even produce a pregear.
\end{prop}

\proof By Proposition \ref{prop:pregearcondition} the condition for
$g(R_0)$ to be a gear is that the image of the upper edge of $\partial
R_0$ should be straight. Parametrize this edge by $\omega_3/2 - s$ for
$s\in[0,\omega_1)$.  Since generally $y_1,y_2$ do not take real values
here, we introduce an auxiliary function
$v(s)=e^{-i\gamma}y(\omega_3/2 - s)$ as the general solution to
\[ 2v'' + \varphi(\frac{\omega_3}{2}-s)v =0. \]
This differential equation has real coefficients, and the initial values are
\[ v(0) = e^{-i\gamma}y(\frac{\omega_3}{2}),\quad
  v'(0) = -e^{-i\gamma}y'(\frac{\omega_3}{2})
\]
where $y$ is the general solution to (\ref{eq:odeR0}).
If we can find $\gamma$ so that both $v(0)$ and $v'(0)$ are real, then
we have $v(t)$ real for all $t$. Clearly the condition is that
$y'(\omega_3/2)$ \textit{is a real multiple of} $y(\omega_3/2)$.
Further, when this holds, we have $y(\omega_3/2-s)=e^{i\gamma}v(s)$,
so $\arg y(\omega_3/2-s)$ is constant.  Therefore $\arg 1/y^2$ is constant on the
top edge of $\partial R_0$.

From (\ref{eq:ggeneral}) we have that $g'=(y_1+\alpha
y_2)^{-2}$.  Writing $y=y_1+\alpha y_2$, our problem is to choose
$\alpha$ so that
\[ y_1'+\alpha y_2' \/\mbox{ \it is a real multiple of }\/  y_1+\alpha y_2 
 \/ \mbox{ \it at } \zeta=\frac{\omega_3}{2}.
\] 
The statement of the Proposition now follows from the fact that given any
$z_1,z_2,w_1,w_2\in\C$, $\alpha\in\R$, then $w_1+\alpha w_2$ is a real
multiple of $z_1+\alpha z_2$ if and only if
\[ (\im z_2 \overline{w_2})\alpha^2 + (\im(z_1 \overline{w_2} + z_2 \overline{w_1}))\alpha  +  \im z_1 \overline{w_1} =0,
\]
as can be seen by multiplying the numerator of $(w_1+\alpha
w_2)/(z_1+\alpha z_2)$ by the conjugate of the denominator and then
taking the imaginary part. 
\qed

\subsubsection{Computational procedure\label{subsubsec:comprect}}
 
Now we can specify the numerical implementation of Proposition
\ref{prop:gearconditiong}.  While (\ref{eq:odeR0}) holds for $\zeta\in
R_0$, numerically we solve it first along the real interval
$[0,\omega_1/2]$ and evaluate the 1-jets of $y_1$, $y_2$ the endpoint
to obtain
\begin{equation}  \label{eq:yright}
 J_{y_1}\left(\frac{\omega_1}{2}\right)=(b_1,b_1'),\quad
  J_{y_2}\left(\frac{\omega_1}{2}\right)=(b_2,b_2').
\end{equation}
Then we solve the initial value problem
 $2u''(s) -\varphi(\omega_1/2 + is) u(s) \ = \ 0$,
 $J_{u_1}(0) = (1,0)$, $J_{u_2}(0) = (0,1)$
for $t\in[0,|\omega_2|/2]$, obtaining at the endpoint
\begin{equation}  \label{eq:uright}
  J_{u_1}\left(\frac{|\omega_2|}{2}\right)= (c_1,c_1'), \quad
  J_{u_2}\left(\frac{|\omega_2|}{2}\right)= (c_2,c_2').
\end{equation}
It follows from this that
$y_1(\omega_1/2+is)=b_1 u_1(s) + ib_1'u_2(s)$,
$y_2(\omega_1/2+is)=b_2 u_1(s) + ib_2'u_2(s)$,
because the left and right sides satisfy identical initial
conditions. Evaluate this at
$s=|\omega_2|/2$ using (\ref{eq:yright}), (\ref{eq:uright}) to obtain
\begin{eqnarray}   
 J_{y_1}\left(\frac{|\omega_3|}{2}\right) &=& (b_1c_1+ib_1'c_2,\ b_1'c_2'-ib_1c_1') ,
  \nonumber\\
 J_{y_2}\left(\frac{|\omega_3|}{2}\right) &=& (b_2c_1+ib_2'c_2,\ b_2'c_2'-ib_2c_1')  .
\label{eq:ycorner} 
\end{eqnarray}
The values $b_1,b_1',b_2,b_2',c_1,c_1',c_2,c_2'$ are all real. They
give in (\ref{eq:ycorner}) the values required in Proposition
\ref{prop:gearconditiong}, which can thus be found by two integrations
on a real interval and two integrations on the right edge of $\partial
R_0$.

\section{Numerical results\label{sec:numerical}} 

We give here some specific results of the numerical application of the
operations described in the previous section.

\subsection{Computing $\beta$ and $\gamma$ as functions of $t$ and $\lambda$}

\noindent\emph{Radial integration.}
The first, most direct method discussed in \ref{subsubsec:diskradial}
above is to substitute into (\ref{eq:2jetfromode}) the boundary values
obtained by solving (\ref{eq:SLt}) along three specific radii in order
to apply Proposition~\ref{prop:curv} and Proposition
\ref{prop:pregeartogear}. Going back and solving
(\ref{eq:SfRtlambda0}) numerically on any radius, we calculate $f(z)$
for any $z\in\D$. In particular, we have all boundary points of the
image gear $f(\D)$ as long as we stay away from the prevertices (where
$R_{t,\lambda}$ has poles).  Then it is easy to find the slope of the tooth
edges, the gear center, and the extreme points
\begin{figure}[!b]   \centering
  \pic{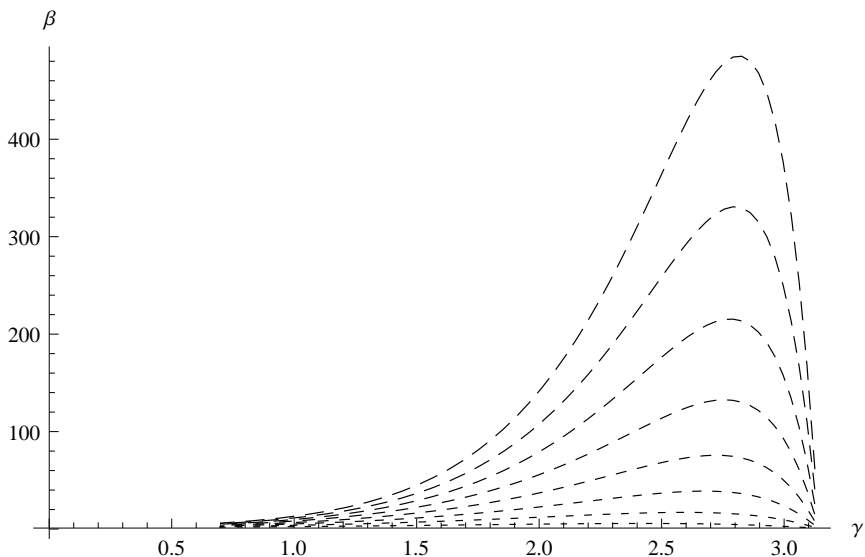}{-1,0}{6.5}{scale=.9}
\caption{Pairs $(\gamma,\beta)$ for fixed $t=\pi/n$ with $n=3,4,\dots,10$.  Larger dashing indicates larger value of $t$.}
\label{fig:gammabetafixedt}
\end{figure}
$f(-1),f(1)$ of the other two edges. From this we have a way of
calculating $\beta,\gamma$ as a function of $t,\lambda$. This method
of calculation was used to produce Figure \ref{fig:gammabetafixedt},
which suggested Conjecture~3.5 of \cite{BrP2} to the effect that for
fixed $\beta$, there are precisely two values of $\gamma$ for each
$t$ below a threshold value $t_\beta$.

\noindent\emph{SPPS integrals.}
The approach described in \ref{subsubsec:diskspps} was carried out
in many examples. To apply Proposition \ref{prop:SPPS}, it is necessary to have a
nonvanishing solution of the differential equation. This amounts to having
a good initial guess for $\lambda_\infty$, which is a simple matter in
the light of formulas (\ref{eq:lambdalimits}) below. 

\begin{figure}[b!]  \centering 
\pic{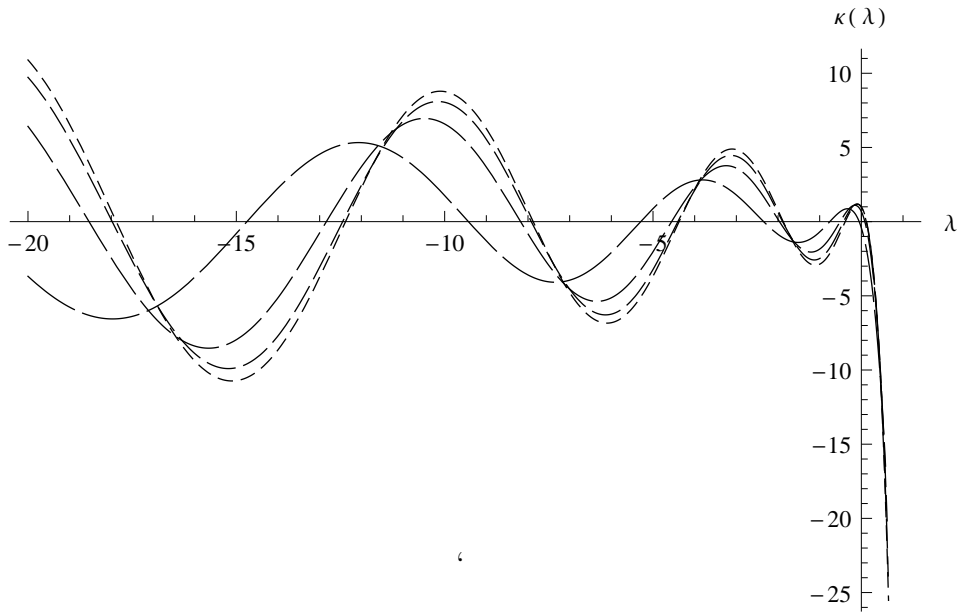}{-1,0}{7}{}`
  \caption{Graphs of the curvature $\kappa(\lambda)$ of the tooth edges of the family of pregears corresponding to $R_{t\lambda}$ for $t/\pi=0.1$, $0.2$, $0,3$, $0,4$.}
  \label{fig:kappagraphs}
\end{figure}

Once this procedure is carried out, we again obtain the data of
Proposition \ref{prop:curv} and Proposition \ref{prop:pregeartogear},
but now expressed as a function of $\lambda$. This permits to study
the behavior of $\beta$ and $\gamma$ as functions of $\lambda$ (for
fixed $t$), and also to solve for $\lambda$ with desired properties.
In particular, Figure \ref{fig:kappagraphs} gives an illustration of
the curvature $\kappa(\lambda)$ of the tooth edges for a fixed value
of $t$. When $\kappa(\lambda)=0$, the mapping $f$ with
$\S_f=R_{t,\lambda}$ sends the arcs $[z_1,z_2]$ and $[z_3,z_4]$ to
subsets of straight lines.  However, in general the solution $f$ is
not univalent: one must choose here the largest value of $\lambda$.
It is interesting that by this method, one can solve for any desired
value (not necessarily zero) and thus obtain pregears with prescribed
radii for the tooth edges.

\begin{figure}[!h]   \centering
 \pic{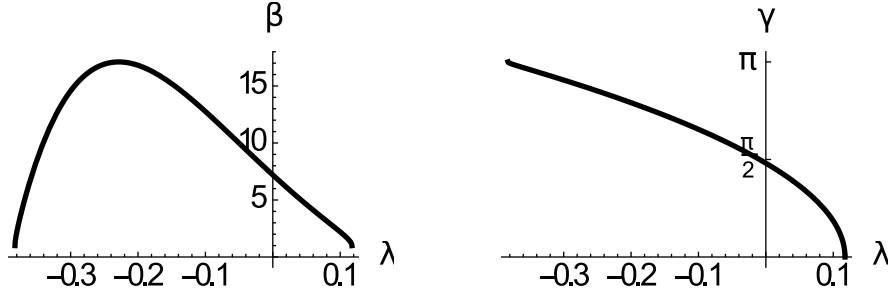}{-1,0}{4}{}  
 \caption{Graphs of  $\beta(\lambda)$ and $\gamma(\lambda)$ produced by SPPS formulas, for $t=\pi/4$. }
 \label{fig:betagamaSPPS}
\end{figure}

As described in \ref{subsubsec:diskspps}, we may substitute the SPPS
formulas into those of Proposition \ref{prop:pregeartogear}, to obtain
the coefficients of the M\"obius transformation $T(z)$ as functions of
$\lambda$.  Applying this ``symbolic'' $T$ to the SPPS series for
$p_{-1}(\lambda)$, $p_1(\lambda)$, $p(\lambda)$, we obtain formulas
$\beta(\lambda)$ and $\gamma(\lambda)$ as combinations of power series
(one could apply the Cauchy rule for products of series, together with
inversion of series, to obtain a single power series in each case, but
this is is rather complicated and fortunately is not necessary). An
example is shown in Figure \ref{fig:betagamaSPPS}.  When combined,
$(\gamma(\lambda),\beta(\lambda))$ parametrize a graph such as in
Figure \ref{fig:gammabetafixedt}.

Inasmuch as this calculation depends on the values of
$R_{t,\lambda}(z)$ for $z$ on the radius from $z=0$ to $z_0=i$ it
necessarily involves complex arithmetic.  As an alternative to this
method, we can also use the other criterion of Proposition
\ref{prop:pregearcondition}: that the non-tooth edges be
concentric. By integrations along the real axis, which only involve
real values of $R_{t,\lambda}$, we can calculate the centers of the
circles containing these edges as $f(-1)+1/\kappa_{-1}$ and
$f(1)+1/\kappa_{1}$, where $\kappa_{-1}$ and $\kappa_{1}$ are the
corresponding curvatures.  Subtracting the SPPS formulas for the
values of the centers produces a function of $\lambda$ which
vanishes when the centers coincide.  Numerical experiments indicate
that either approach seems to work equally well.

\noindent\emph{Repositioning of gear center.}
To avoid possible confusion we reiterate that with the condition
$J_f(0)=(0,1,0)$, the value $f(0)$ is not likely to be the gear center
$w_0$ of $f(\D)$. Therefore such mappings are quite different from the
solutions $F$ of the classical integral representation of gear
mappings as in \cite{BPea,Goo,Pea}.  Consider a general configuration
of prevertices $\pm e^{\pm it_1},\pm e^{\pm it_2}$ (recall the
discussion at the beginning of \ref{subsec:diskschw}), and let $f$ be
the solution of $\S_f=R_{t_1,t_2,\lambda}$ normalized by
$J_f(0)=(0,1,0)$.  Having calculated the data of Proposition
\ref{prop:pregeartogear} we can find $w_0$, and presumably having
already solved $f=y_2/y_1$ along $[-1,1]$, it is possible now to
approximate $p=f^{-1}(w_0)$ by numerical inversion.  Then the
composition $F=f\circ T_{-p}$ (which has a different 2-jet at the
origin) sends 0 to $w_0$.  Further, $F(\D)=f(\D)$.  Clearly
$S_F=R_{t_1',t_2',\lambda}$ where $T_{-p}(e^{it_1'})=e^{it_1}$,
$T_{-p}(e^{it_2'})=e^{it_2}$.
Based on numerical examples of this procedure we are led to conjecture 
the following.

\begin{conj} \label{conj:uniquelambda} Given $t_1,t_2$,
  $0<t_1<t_2<\pi/2$, there is a unique $\lambda\in\R$ such that the
  solution $f$ of $\S_f=R_{t_1,t_2,\lambda}$ normalized by
  $J_f(0)=(0,1,0)$ is a gear mapping.
\end{conj}
This is equivalent to the statement that if there is a conformal
mapping of gears $G_{\beta,\gamma}\to G_{\beta',\gamma'}$ respecting
the vertices and with 2-jet of the form $(0,r,0)$ at the origin
($r>0$), then $\beta=\beta'$ and $\gamma=\gamma'$.

\begin{figure}[h!tb]   \centering
 \pic{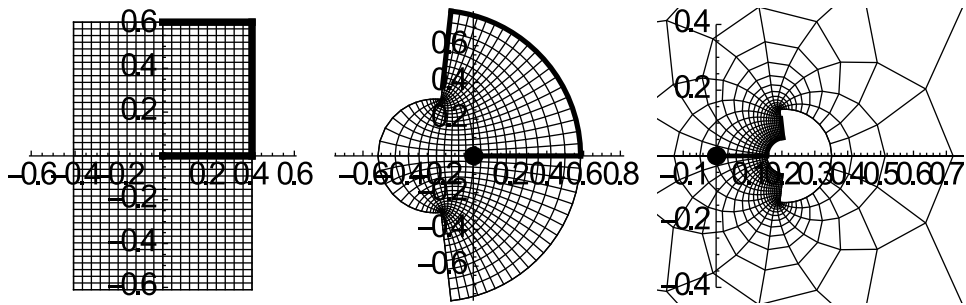}{-1.5,0}{4}{scale=1}  
 \caption{  Mapping of rectangle to gear ($\tau=1.5$).   }
 \label{fig:rectanglemap}
\end{figure}

\noindent\emph{Integration in a rectangle.} The method described in
\ref{subsubsec:comprect} is easily applied, integrating from 0 to
$\omega_1/2$ and from $\omega_1/2$ to $(\omega_1+\omega_2)/2$ to find
the required parameter $\alpha$ for converting the gear to a
pregear. In the first integration one may save the values of $y_1,y_2$
at points along $[0,\omega_1]$ and then use them as initial values for
integrating upwards or downward on vertical segments passing through
those points.  An example is shown in Figure \ref{fig:rectanglemap}.
The image at the right results from the value of $\alpha$ which
produces an unbounded domain. The complement of this unbounded domain
is again a gear domain; we know of no relation between this and the
bounded image.

 \subsection{Inverse problem: prescribed gear parameters}\label{inverseproblem}

 As occurs in many contexts in conformal mapping, the interesting
problem is to find the auxiliary parameters which produce a given
geometry.  First we consider the following simpler question: Given
$t$, determine $\lambda$ so that the image of the gear mapping
$f_{t,\lambda}$ has gear ratio $\beta$ (or alternatively, gear angle
$\gamma$).  Since $t$ is fixed, we can calculate the formula
approximating $\kappa(\lambda)$, solve $\kappa(\lambda) = 0$, and then
having $t,\lambda$ we determine $\beta$.  This can be repeated as
necessary, by a process of successive approximations, to make $\beta$
have the desired value within specified accuracy. We have carried out
this approach successfully.  A more direct method is to use the SPPS
formulas for $\beta(\lambda)$, $\gamma(\lambda)$ (recall Figure
\ref{fig:betagamaSPPS}), and solve them directly for the desired
values of $\beta,\gamma$. However, for extremely small values of $t$
this does not give good results unless a great number of powers are
taken in the SPPS formulas.

\begin{figure}[!b]   \centering
 \pic{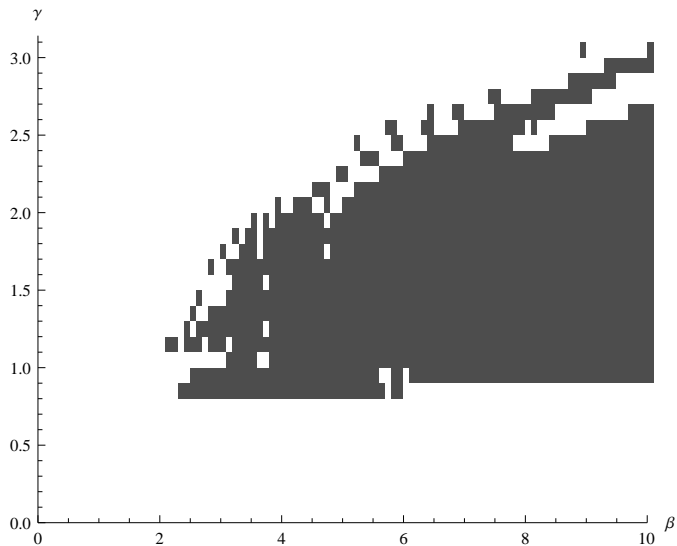}{0,0}{7}{scale=.7}  
 \caption{Experimental region of success of Broyden's method for inverting the
   parameter correspondence, applied for a rectangular grid of
   $1.1\le\beta\le10.0$ and $0.1\le\gamma\le\pi-0.1$.  }
 \label{fig:failbroyden}
\end{figure}

The method of Broyden \cite{FKS} may also be used to find zeros of the
mapping $(t,\lambda)\to(\beta-\beta_0,\gamma-\gamma_0)$, without
recourse to the SPPS formulas. This turns out to be extremely fast
(thousands of solutions in less than a second on an ordinary portable
computer).  However, the Broyden method works by jumping around
unpredictably in $\R^2$ from the initial guess for $(t,\lambda)$, and
may fail by leaving the $(t,\lambda)$ region where $(\beta,\gamma)$ is
well defined.  Figure \ref{fig:failbroyden} shows values where this
method succeeds starting from the initial guess
$(t,\lambda)=(\pi/4,0)$.  We will not pursue further the question
of improving the initial guess.

\section{Internal structure of the region of gearlikeness}\label{sec:internal}

In \cite{BrP2} it was shown that the region of gearlikeness in the
$(t,\lambda)$-plane is
\[ \G = \{ (t,\lambda)\colon\ \lambda_t^-<\lambda<\lambda_t^+ \} ,
\] where
\begin{eqnarray}  \label{eq:lambdalimits}
 \lambda_t^- = -\frac{1}{4} - \frac{1}{16}\left(\cos t +
    \frac{1}{\cos t}\right) ,\quad
     \lambda_t^+ = \frac{1}{4} - \frac{1}{16}\left(\cos t + 
    \frac{1}{\cos t}\right).
\end{eqnarray}
  
We use our numerical methods to obtain a very
illuminating picture of the structure of $\G$ as related to
$(\beta,\gamma)$.  Observe that the vertical
cross-sections of $\G$ (representing gears of a given conformal
modulus $M(t)$) are of common height $\lambda_t^+ - \lambda_t^-=1/2$
for all $t$, and
\[ \lambda_0^-=\lim_{t\to0} \lambda_t^- = -\frac{3}{8}; \quad
   \lambda_0^+=\lim_{t\to0} \lambda_t^+ = \frac{1}{8}.
\]   

\begin{figure}[!tb]  \centering
  \pic{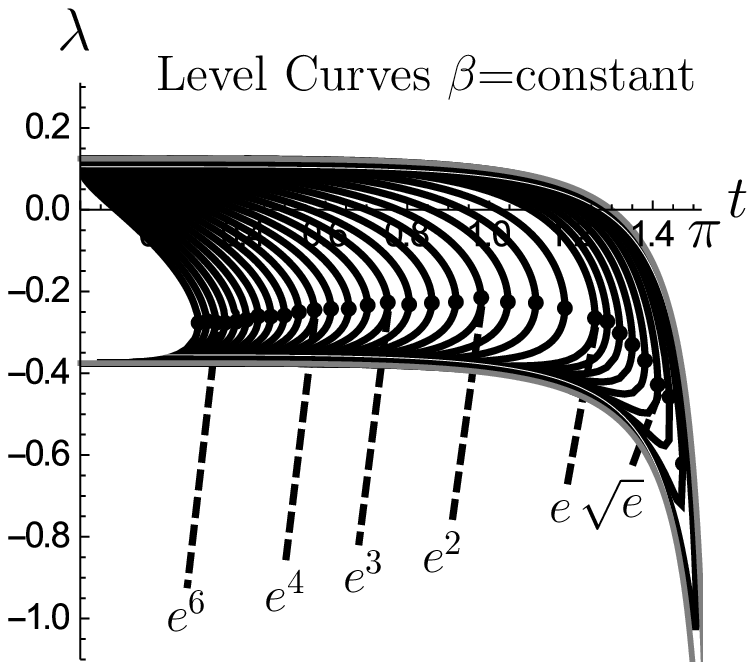}{-1.5,0}{4}{scale=.8}
  \pic{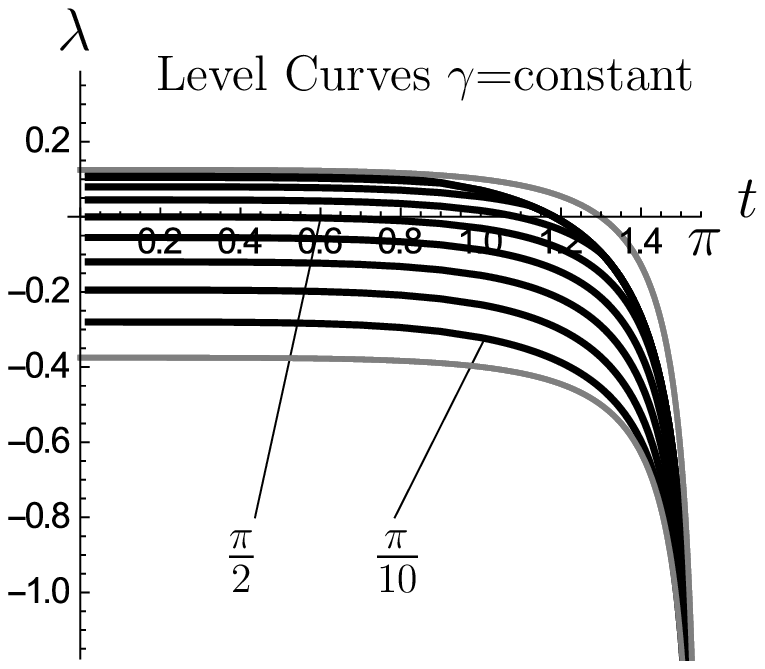}{5.7,0}{4}{scale=.8}
  \caption{Region of gearlikeness $\G$ foliated by $\beta$ level curves
(left) and $\gamma$ level curves (right).  }
  \label{fig:gearregionstructure}
\end{figure}

We use the methods of the previous section for $(\beta,\gamma)\mapsto(t,\lambda)$
to calculate the level curves of the geometric parameters.
Figure \ref{fig:gearregionstructure} shows the subsets of $\G$ of
constant $\log\beta=$ 0.2, 0.4, 0.6, 0.8, 1.0; 1.25, 1.5, 2.0\dots\ and of constant $\gamma=0.1\pi,\ 0.2\pi,\ \dots, 0.9\pi$.  

Individually, the $\beta$-curves accumulate only at the two extreme
points $(0,\lambda_0^-)$ and $(0,\lambda_0^+)$, never at interior
points $(t,\lambda_t^\pm)$ of the lower and upper boundaries.  At
these extreme points $\beta\to1$ while $\gamma\to\pi,0$ respectively,
and the gear $G_{\beta,\gamma}$ degenerates to a disk in either case.
It may appear paradoxical that the limiting Schwarzian derivatives
$R_{t,\lambda_t^\pm}(z)$ are not identically zero; however, the
pullbacks of the Schwarzian derivatives according to Proposition
\ref{prop:zerotozero} by appropriate $T_q$ (with $q$ depending on
$\lambda$) do vanish in the limit.  For fixed $\gamma$, as $t\to0$ or
$t\to\pi/2$ we have $\beta\to\infty$ or $\beta\to1$ respectively.  For
fixed $\beta$, as $t\to0$ we have $\lambda\to\lambda_0^-$ or
$\lambda\to\lambda_0^+$ and then $\gamma\to\pi$ or
$\gamma\to0$.  This is one of the qualitative results on the
conformal module proved in \cite{BrP2}.
 
Each level curve for fixed $\gamma$ intersects the $\lambda$-axis in a
value, which might be termed $\lim_{t\to0}\lambda^{[t,\gamma]}$,
corresponding to a degenerate gear. As $\beta\to\infty$ for fixed $\gamma$, the gear $G_{\beta,\gamma}$
tends to the union of $\D$ with a full sector of angle
$2\gamma$. Letting $t_1\to0$ in the gear mapping, we are led to
apply the general formula (4) to this circular triangle with angles
$3\pi/2$, $3\pi/2$, $\gamma$ and prevertices $e^{\pm i t_2}$, $1$, to
obtain
\begin{eqnarray*}
  \S_f(z) &=& \frac{(2\gamma/\pi)^2(\cos t_2-1)}{(z-1)^2(z^2-(2\cos t_2)z+1)}\\
    &&\ \  -\  \frac{ (\cos t_2-1)((5 z^2-14 z+5)\cos t_2+7z^2-10z+7)} 
          {2(z-1)^2(z^2-(2\cos t_2)z+1)^2}
\end{eqnarray*}  
whereas
\begin{eqnarray*}
 R_{0,t_2,\lambda}(z) =  \frac{8\lambda(\cos t_2-1)}{(z-1)^2(z^2-(2\cos t_2)z+1)} -
     \frac{5\sin^2t_2} {2(z-1)^2(z^2-(2\cos t_2)z+1)^2}.
\end{eqnarray*}
Equating these two Schwarzian derivatives we find that
\begin{equation} \label{eq:limitlambda}
 \lim_{t\to0}\lambda^{[t,\gamma]}   =
    \frac{1}{8}\left( 1 - \left(\frac{2\gamma}{\pi}\right)^2\right)
\end{equation}
independently of $t_2$. This formula is confirmed by Figure
\ref{fig:gearregionstructure} (right).

Goodman noted in \cite{Goo} that when one limits the discussion to a
particular $\gamma$ (and fixes the normalization as $f(0)=0=f'(0)-1$),
a relation is determined between the parameters $t_1$ and $t_2$,
stating that he could not calculate it except for the particular case
$\gamma=\pi/2$, where he found that
\[ \cos t_1=1,\quad \cos t_2=\frac{1}{2},
\]
i.e.\ $t_1=0$, $\beta=\infty$ as the outer vertices of the gear have
coalesced at $\infty$, while $t_2=\pi/3$.  Using Goodman's explicit
mapping formula
\[  f(z) = \frac{4}{27}\frac{2(1-z+z^2)-2+3z+3z^2-2z^3}{z(1-z)}
\] and comparing $\S_f$ with $R_{0,\pi/3,\lambda}$, one finds readily
that $\lambda=0$, which is thus the value of  $\lim_{t\to0}\lambda^{[t,\pi/2]}$,
confirming  (\ref{eq:limitlambda}) for this case.

\section{Applications and conclusions}\label{sec:appl}

We close with some brief applications of our results on gear mappings.

\subsection{The first Maclaurin series coefficient}\label{subsec:maclaurin}

In \cite{Goo} Goodman left unsolved the problem of calculating the
ratio $b_1/f(1)$, where $b_1=f'(0)$ is the first Maclaurin coefficient of the conformal
mapping $h$ of $\D$ onto a one-tooth gear domain that maps the origin to the
gear center. This ratio was expressed in terms of singular
integrals in \cite{Pea} and later in \cite{Br3}, where it is shown
that 
\[ f'(0)=f(1)\int_0^1 \left(  \frac{1}{x} - 
 \frac{\sqrt{1-(\cos t_2)x}}{x\sqrt{1-(\cos t_1)x}\sqrt{1-x^2}} \right)\,dx  .
\]
The integrand presents rather complicated singularities at the endoints
of integration.
 
In the construction of Proposition \ref{prop:zerotozero} we obtained
the 2-jet $J_f(0)$ via (\ref{eq:2jetzerotozero}), so we have
$f'(0)$. It is also a simple matter to calculate $f(1)$
numerically. Thus one may evaluate integrals of the form given above
via solutions of gear mapping problems.


\subsection{Module of the complement of annular rectangle}

\begin{figure}[!tb]  \centering
  \pic{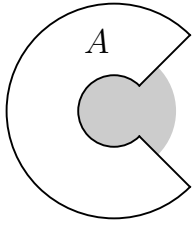}{4,0}{2}{scale=.9}
 \caption{Complement of annular rectangle is formed of a gear and
its reflection along its B-arc. }
  \label{fig:annularrectangle}
\end{figure}

Consider the bounded region $A=A_{\beta,\gamma}$ with boundary
$ \partial A = \{e^{i\theta}\colon\ \gamma<\theta<2\pi-\gamma\}\cup
\{re^{i\gamma}\colon\ 1\leq r \leq \beta^2\} \cup
\{\beta^2e^{i\theta}\colon\ \gamma<\theta<2\pi-\gamma\}
\cup\{re^{-i\gamma}\colon\ 1\leq r\leq \beta^2\}, $ i.e., an ``annular
rectangle'' as in Figure \ref{fig:annularrectangle} subtending an
angle $2(\pi-\gamma)$ within an annulus of radii of ratio $\beta^2$.
While it is easy to map a disk or rectangle to $A$, it is not known
(cf.\ \cite[p.\ 122]{Ku}) how to obtain an expression for the
conformal module of the exterior $A^*$ of $A$ in closed form in terms
of of the modulus of $A$. (This is the situation with quadrilaterals
in general. There has been a surge of interest recently \cite{DP,
  HRV1, HRV2, VZ} in the question of numerical calculation of exterior
moduli of topological quadrilaterals.) Here we give a numerical
solution to this problem.

\begin{theo}\label{thm:annularrectangle}
  Let $A$ be the annular rectangle with inner radius $1$ and outer
  radius $\beta^2$ and angle $2(\pi-\gamma)$. The conformal module of
  the exterior $A^*=(\C\cup\{\infty\})\setminus \mbox{cl}\,A$ is half the
  module of the gear domain with gear angle $\gamma$ and gear ratio
  $\beta$.
\end{theo}
\proof Consider the conformal mapping $g\colon R_0\to
G_{\gamma,\beta}$ of Proposition~\ref{prop:recgear}. The Schwarz
reflection applied across the right vertical edge of $R_0$ produces a
mapping onto $A^*$ from a rectangle having double the width of $R_0$.
\qed

Given an annular rectangle $A_{\beta,\gamma}$, the values $t$ and
$\lambda$ corresponding to $\beta$ and $\gamma$ can be computed as
described in Section~\ref{inverseproblem}. This means that
$f_{t,\lambda}$ is a mapping onto a gear $G_{\beta,\gamma}$ with
conformal module $M(t)$. According to
Theorem~\ref{thm:annularrectangle}, the module of $A_{\beta,\gamma}$
is therefore calculated numerically to be $M(t)/2$ .
 
\subsection{Multitooth gears} 

 Let $f\colon\D\to G_{\beta,\gamma}$ be a conformal mapping to a normalized gear
domain, $f(0)=0$.  The function $f_n(z) = \sqrt[n]{f(z^n)}$ is a
mapping to a regular $n$-toothed gear. Applying the Chain Rule to
$P_n\circ f_n = f \circ P_n$ where $P_n(z)=z^n$, and using $\S_{P_n}=
(1-n^2)z^{-2}/2$ we find
 \[  \frac{1-n^2}{2}\frac{f_n'(z)^2}{f_n(z)^2} +\S_{f_n}(z) =   
  \S_f(z^n)\, n^2z^{2(n-1)} + \frac{1-n^2}{2}{z^2}.
\]
Since both $\S_f$ and $\S_{f_n}$ are rational functions, it follows
that $(f_n'/f_n)^2$ is a rational function. With somewhat more work
one recovers the formula of \cite{Goo} of which equation (2) of \cite{BrP2} is
a particular case.

\begin{figure}[!t]  \centering
  \pic{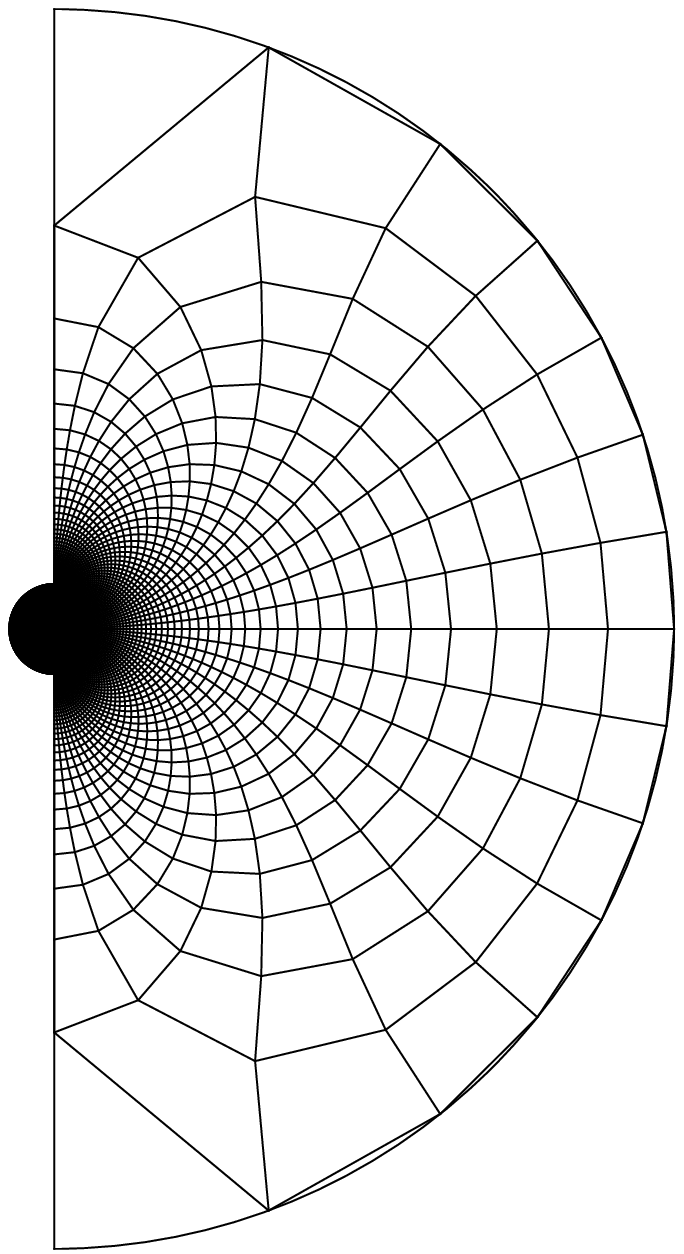}{-2,0}{4.5}{scale=.4}
  \pic{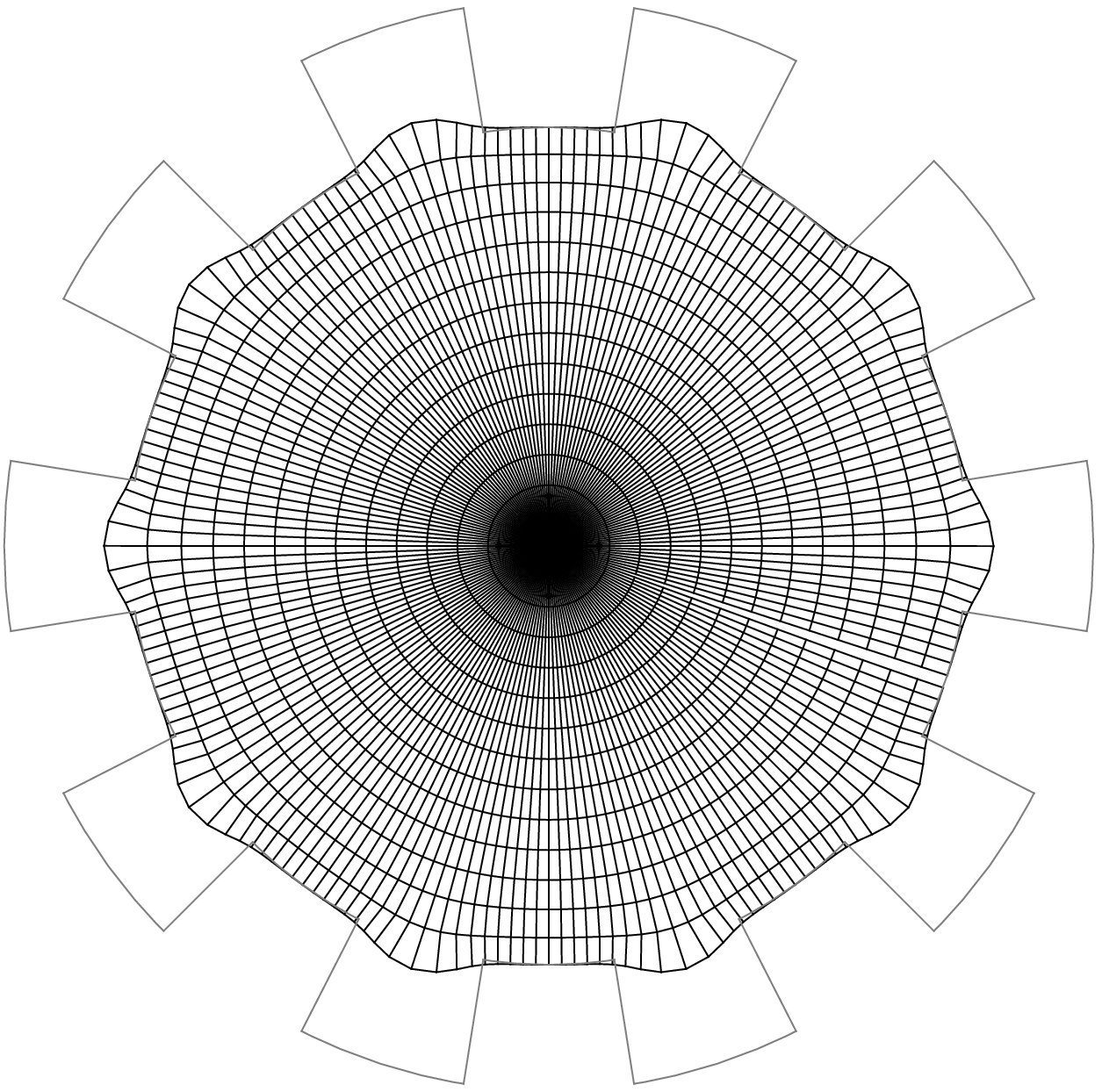}{1.8,-.3}{4}{scale=.45}
  \pic{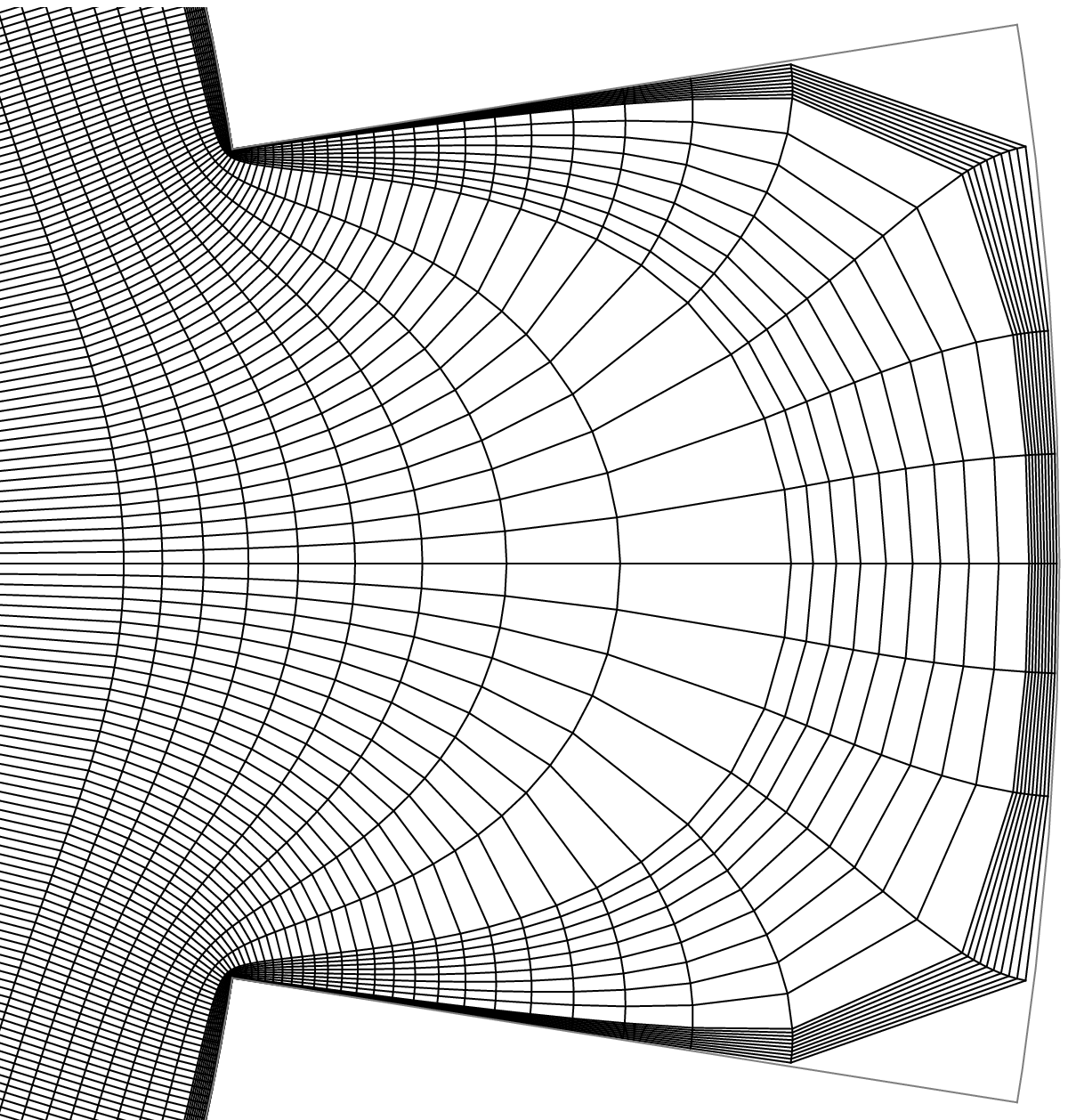}{7.5,.4}{4}{scale=.3}
  \caption{Single-tooth gear (left) calculated to generate a multitooth
    version of prescribed geometry (center). Detail (right) includes level curves corresponding to
    $|z|=0.9,\ 0.91,\dots,\ 0.99,\ 0.991,\dots,.999,\
    0.9991,\dots,0.9999$.}
  \label{fig:multitooth}
\end{figure}

The image $f_n(G_{\beta,\gamma})$ has gear ratio $\beta^{1/n}$, and
each tooth subtends an angle of $2\gamma/n$ with spacing of
$2(\pi-\gamma)/n$ between consecutive teeth. By means of these facts
it is simple to calculate any desired regular multitooth domain.  We
illustrate with the example of $n=10$ teeth, having intertooth space
equal to the tooth width, and gear ratio arbitrarily chosen as
1.3. Thus we have $\beta=1.3^{10}\approx13.79$,
$\gamma\approx\pi/2$. By the methods given in Section~\ref{inverseproblem}, we find that the one-tooth gear with these
parameters is obtained by $t\approx0.6024$,
$\lambda\approx-0.0029$. The result is in Figure
\ref{fig:multitooth}. We stress that the renormalization worked out in
Proposition~\ref{prop:zerotozero} is essential, since without the
condition that $f(0)=0$ is the gear center, the application of $P_n$
will not work.
 
It may be noted that when the number $n$ of teeth is large, the circular arcs
are approximated by straight lines, and the gear mapping may be approximated by
a Schwarz=Christoffel integral.

\appendix
\section{Appendix: SPPS method} 
The sequence $I_n$ of \textit{iterated integrals} generated by an
arbitrary pair of functions $(q_0,q_1)$ is defined recursively by
setting $I_0=1$ identically and for $n\ge1$,
\begin{equation}
  I_n(z) = \int_0^z I_{n-1}(\zeta)\,q_{n-1}(\zeta)\,d\zeta
\end{equation}
where $q_{n+2j}=q_n$ for $j=1,2,\dots$

\begin{prop}{\rm \cite{KP1}} \label{prop:SPPS}
 Let\/ $\psi_0$ and $\psi_1$ be given, and suppose that\/ $y_\infty$
is a nonvanishing solution of
\[ y_\infty''+\psi_0\,y_\infty=\lambda_\infty\psi_1\,y_\infty
\] on the interval $[0,1]$, where $\lambda_\infty$ is any constant.
Choose\/ $q_0=1/y_\infty^2$, $q_1=\psi_1\,y_\infty^2$ and define\/
$\X{n}$, $\Xt{n}$ to be the two sequences of iterated integrals
generated by\/ $(q_0,q_1)$ and by\/ $(q_1,q_0)$, respectively.  Then
for each\/ $\lambda\in\C$ the functions
\begin{eqnarray}
  y_1 &=& y_\infty \sum_{k=0}^\infty(\lambda-\lambda_\infty)^k
            \Xt{2k}, \nonumber\\ 
  y_2 &=& y_\infty \sum_{k=0}^\infty(\lambda-\lambda_\infty)^k
              \X{2k+1} \label{eq:SPPSseries}
\end{eqnarray}
are linearly independent solutions of the equation
\begin{equation}  \label{eq:y''}
 y''+\psi_0y=\lambda\psi_1y 
\end{equation}
on  $[0,1]$.  Further, the series for $y_1$ and $y_2$ 
converge uniformly on $[0,1]$ for every $\lambda$.
\end{prop}
 
It is a straightforward matter to obtain the appropriate linear
combination of solutions with desired 2-jets at $z=0$, for example
$(0,1)$ and $(1,0)$.

\noindent Philip R. Brown\\ 
Department of General Academics\\
Texas A\&M University at Galveston\\ 
PO Box 1675, Galveston, Texas 77553 -1675 \\ 
\texttt{ brownp@tamug.edu }

\medskip
\noindent R. Michael Porter \\
Departamento de Matem\'aticas, CINVESTAV--I.P.N.\\
Apdo.\ Postal 1-798, Arteaga 5 \\
Santiago de Queretaro, Qro., 76000  MEXICO \\
\texttt{mike@math.cinvestav.edu.mx}


\end{document}